\documentclass[a4paper,12pt]{amsart}
\usepackage{amssymb}
\usepackage{ifthen}
\usepackage{graphicx}
\nonstopmode \numberwithin{equation}{section}
\newtheorem{thm}[equation]{Theorem}
\newtheorem{cor}[equation]{Corollary}
\newtheorem{lem}[equation]{Lemma}

\newtheorem{rem}[equation]{Remark}
\newtheorem{rems}[equation]{Remarks}
\theoremstyle{definition}
\newtheorem{defin}[equation]{Definition}
\newtheorem{examp}[equation]{Example}
\newtheorem{prob}[equation]{Problem}
\newtheorem{ques}[equation]{Question}

\newtheorem{conj}[equation]{Conjecture}
\newtheorem{deter}[equation]{Determination}
\newtheorem{case}[equation]{Case}
\newtheorem{claim}[equation]{Claim}

\newcounter {own}
\def\theown {\thesection       .\arabic{own}}

\newenvironment{pf}[1][]{%
 \vskip 3mm
 \noindent
 \ifthenelse{\equal{#1}{}}%
  {{\bf Proof. }}%
  {{\bf #1.} }%
 }%
{\qed\bigskip}

\newcounter{alphabet}
\newcounter{tmp}

\newcommand{\F}{{\mathcal{F}}}

\newcommand{\IC}{{\mathbb C}}

\newcommand{\M}{{\mathrm M}}

\def\be{\begin{equation}}
\def\ee{\end{equation}}

\newcommand{\bee}{\begin{enumerate}}
\newcommand{\eee}{\end{enumerate}}

\newcommand{\blem}{\begin{lem}}
\newcommand{\elem}{\end{lem}}
\newcommand{\bthm}{\begin{thm}}
\newcommand{\ethm}{\end{thm}}
\newcommand{\bcor}{\begin{cor}}
\newcommand{\ecor}{\end{cor}}
\newcommand{\beg}{\begin{examp}}
\newcommand{\eeg}{\end{examp}}
\newcommand{\begs}{\begin{examples}}
\newcommand{\eegs}{\end{examples}}
\newcommand{\bdefe}{\begin{defin}}
\newcommand{\edefe}{\end{defin}}
\newcommand{\bprob}{\begin{prob}}
\newcommand{\eprob}{\end{prob}}
\newcommand{\bques}{\begin{ques}}
\newcommand{\eques}{\end{ques}}
\newcommand{\bei}{\begin{itemize}}
\newcommand{\eei}{\end{itemize}}

\newcommand{\bde}{\begin{deter}}
\newcommand{\ede}{\end{deter}}
\newcommand{\bca}{\begin{case}}
\newcommand{\eca}{\end{case}}
\newcommand{\bcl}{\begin{claim}}
\newcommand{\ecl}{\end{claim}}
\newcommand{\bcon}{\begin{conj}}
\newcommand{\econ}{\end{conj}}
\newcommand{\bcons}{\begin{conjs}}
\newcommand{\econs}{\end{conjs}}
\newcommand{\bprop}{\begin{propo}}
\newcommand{\eprop}{\end{propo}}
\newcommand{\br}{\begin{rem}}
\newcommand{\er}{\end{rem}}
\newcommand{\brs}{\begin{rems}}
\newcommand{\ers}{\end{rems}}
\newcommand{\bo}{\begin{obser}}
\newcommand{\eo}{\end{obser}}
\newcommand{\bos}{\begin{obsers}}
\newcommand{\eos}{\end{obsers}}
\newcommand{\bpf}{\begin{pf}}
\newcommand{\epf}{\end{pf}}
\newcommand{\ba}{\begin{array}}
\newcommand{\ea}{\end{array}}
\newcommand{\beq}{\begin{eqnarray}}
\newcommand{\beqq}{\begin{eqnarray*}}
\newcommand{\eeq}{\end{eqnarray}}
\newcommand{\eeqq}{\end{eqnarray*}}

\title{Mapping problems for quasiregular mappings}

\author{Manzi Huang}
\address{Manzi Huang, Department of Mathematics,
Hunan Normal University, Changsha,  Hunan 410081, People's Republic
of China}
\email{mzhuang79@yahoo.com.cn}

\author{Antti Rasila$^*$}
\address{Antti Rasila,
 Department of Mathematics, Hunan Normal University, Changsha,  Hunan 410081, People's Republic
of China; Department of Mathematics and Systems Analysis,
Aalto University, FI-00076 Aalto, Finland}
\email{antti.rasila@iki.fi}

\author{Xiantao Wang}
\address{Xiantao Wang, Department of Mathematics,
Hunan Normal University, Changsha,  Hunan 410081, People's Republic
of China} 
\email{xtwang@hunnu.edu.cn}

\subjclass[2000]{Primary:
30C65, 30F45; Secondary: 30C20} \keywords{Quasiregular mapping,
closed quasiregular mapping, property $P_1$, property $P_2$, maximal
(minimal) multiplicity.}

\thanks{$^*$The second author is the
corresponding author}

\begin{document}

\begin{abstract}
We study images of the unit ball under certain special classes of quasiregular mappings. For homeomorphic, i.e., quasiconformal mappings problems of this type have been studied extensively in the literature. In this paper we also consider non-homeomorphic quasiregular mappings. In particular, we study (topologically) closed quasiregular mappings originating from the work of J. V\"ais\"al\"a and M. Vuorinen in 1970's. Such mappings need not be one-to-one but they still share many properties of quasiconformal mappings. The global behavior of closed quasiregular mappings is similar to the local behavior of quasiregular mappings restricted to a so-called normal domain. 
\end{abstract}

\maketitle

\section{Introduction}

We consider quasiregular mappings in the $n$-dimensional Euclidean space $\mathbb{R}^n$. Quasiconformal and quasiregular mappings in $\mathbb{R}^n$, $n\ge 3$ are natural generalizations of conformal and analytic functions of one complex variable, respectively. For basic properties of these classes of mappings, we refer to \cite{Rickman:1993,Vaisala:1971,Vuorinen:1988}. In the complex plane, it follows from the Riemann mapping theorem that any simply connected domain is the image of the unit disk in a conformal, and thus quasiconformal, mapping. The so-called measurable Riemann mapping theorem further generalizes this result by allowing one to find a quasiconformal mapping of given dilatation. However, the problem of characterizing the quasidisks, i.e., quasiconformal images of the unit disk in the quasiconformal mappings of the the whole plane onto itself is interesting (see e.g. \cite{Beurling-Ahlfors,Gehring-Hag}). For $n\ge 3$ the question of characterizing the quasiconformal images of the unit ball $\mathbb{B}^n$ is highly non-trivial, and it has been studied by many authors \cite{Gehring-Vaisala,Hag-Vamanamurthy,Vaisala:1961}. In this paper, we present several examples related to this topic, and new results concerning the so-called closed quasiregular mappings.

The topological properties of quasiregular mappings are similar to those of analytic functions. It is well-known that a nonconstant quasiregular mapping is discrete (i.e. sets $f^{-1}(y)$ are discrete) and open (see e.g. \cite[I.4.1]{Rickman:1993}). We study a subclass of the quasiregular mappings which are characterized by the property that they preserve closed sets. This class of mappings is more general than the quasiconformal mappings, as closed mappings need not be homeomorphic. The class of closed quasiregular mappings originates from the work of J. V\"ais\"al\"a \cite{Vaisala:1966} and M. Vuorinen \cite{Vuorinen:1976,Vuorinen:1978,Vuorinen:1985}.

The global behavior of closed quasiregular mappings is similar to the behavior of quasiregular mappings restricted to the so-called normal domains. The existence of such neighborhoods is well-known, but usually nothing is known of their diameter. The importance of the assumption that mappings are closed arises from the fact that it allows us to extend local estimates which are based on the conformal modulus to global ones.

\section{Preliminaries}

We shall follow standard notation and terminology adopted from
 \cite{Vaisala:1971}, \cite{Vuorinen:1988} and \cite{Rickman:1993}. For
 $x\in\mathbb{R}^n$, $n\geq2$, and $r>0$, let $\mathbb{B}^n(x,r)=\{z\in
 \mathbb{R}^n:|z-x|<r\}$,
$\mathbb{S}^{n-1}(x,r)=\partial\mathbb{B}^n(x,r)$,
$\mathbb{B}^n(r)=\mathbb{B}^n(0,r)$,
$\mathbb{S}^{n-1}(r)=\partial\mathbb{B}^n(r)$,
$\mathbb{B}^n=\mathbb{B}^n(1)$ and
$\mathbb{S}^{n-1}=\partial\mathbb{B}^n$. The space
$\overline{\mathbb{R}}^n=\mathbb{R}^n\cup\{\infty\}$ is the
one-point compactification of $\mathbb{R}^n$. The surface area of
$\mathbb{S}^{n-1}$ is denoted by $\omega_{n-1}$ and $\Omega_n$ is
the volume of $\mathbb{B}^n$. It is well-known that
$\omega_{n-1}=n\Omega_n$ and that
$$
\Omega_n=\frac{\pi^{n/2}}{\Gamma(1+n/2)}
$$
for $n=2,3,\ldots$, where $\Gamma$ is Euler's gamma function. The
standard coordinate unit vectors are denoted by $e_1,\ldots,e_n$.
The Lebesgue measure on $\mathbb{R}^n$ is denoted by $m$.

\subsection*{Quasiregular mappings}

A continuous mapping $f\colon G\to\mathbb{R}^n$, $n\geq2$, of a domain $G$ in
$\mathbb{R}^n$ is called {\it quasiregular} if $f$ is in the Sobolev
space $W^{1,n}_{\mathrm{loc}}(G)$, and there
 exists a constant $K$, $1\leq K<\infty$, such that the inequality
$$
|f'(x)|^n\leq KJ_f(x)
$$
holds a.e. in $G$, where $f'(x)$ is the formal derivative of $f$, and
$|f'(x)|=\max_{|h|=1}|f'(x)h|$. The smallest $K\geq 1$ for which
this inequality is true is called the outer dilatation of $f$ and
denoted by $K_O(f)$. If $f$ is quasiregular, then the smallest
$K\geq1$ for which the inequality
$$
J_f(x) \leq K l(f'(x))^n
$$
holds a.e. in $G$ is called the inner dilatation of $f$ and denoted
by $K_I(f)$, where $l(f'(x))=\min_{|h|=1}|f'(x)h|$. The maximal
dilatation of $f$ is the number $K(f)=\max\{K_I(f),K_O(f)\}$. If
$K(f)\leq K$, $f$ is said to be $K$-quasiregular. A quasiregular
homeomorphism $f\colon G\to fG$ is
 called
{\it quasiconformal}.

By generalized Liouville's theorem for $n \geq 3$, every $1$-quasiregular
mapping in $\mathbb{R}^n$ is a restriction of a M\"obius
transformation or a constant.
 The M\"obius transformations are very useful in the study of
quasiregular mappings. In particular, we make use of the mapping
$T_a$, $a\in\mathbb{B}^n$, which is the M\"obius transformation with
$T_a(\mathbb{B}^n)=\mathbb{B}^n$, $T_a(a)=0$ and for $e_a=a/|a|$,
$T_a(e_a)=e_a$ and $T_a(-e_a)=-e_a$. For $a=0$, we set
$T_0=\mathrm{id}$ (see \cite[p. 11]{Vuorinen:1988} or \cite[II
2.6]{Ahlfors:1981}).

\subsection*{Modulus of a path family}

Let $\Gamma$ be a path family in $\mathbb{R}^n$, $n\geq 2$. Let
$\F(\Gamma)$ be the set of all Borel functions $\rho\colon
\mathbb{R}^n\to[0,\infty]$ such that
$$
\int_\gamma \rho\,ds\geq 1
$$
for every locally rectifiable path $\gamma\in \Gamma$. The functions
in $\F(\Gamma)$ are called {\it admissible} for $\Gamma$. For $1\leq
n < \infty$, we define
\begin{equation}
\label{modulus} \M(\Gamma) = \inf_{\rho \in
\F(\Gamma)}\int_{\mathbb{R}^n}\rho^n\,dm
\end{equation}
and call $\M(\Gamma)$ the {\it $($conformal$)$ modulus} of $\Gamma$.
If $\F(\Gamma)=\emptyset$, which is true only if $\Gamma$ contains
constant paths, we set $\M(\Gamma)=\infty$. If $\Gamma_1,\Gamma_2$ are path
families in $\mathbb{R}^n$, and every $\gamma\in \Gamma_2$ has a subcurve
in $\Gamma_1$, we say that $\Gamma_2$ is {\it minorized} by $\Gamma_1$ and write
$\Gamma_2>\Gamma_1$. If $\Gamma_1 < \Gamma_2$, then $\M(\Gamma_1)\geq \M(\Gamma_2)$.
For the basic properties of the modulus of the path family, we refer to \cite{Rickman:1993,Vaisala:1971,Vuorinen:1988}. It is well-known that the modulus of a path family is invariant
under conformal mappings. We denote by $\Delta(A,B; G)$ the family of
paths joining $A$ and $B$ in $G$.

We use the following well-known identity of the modulus of the
 spherical annulus: Let $0<a<b$. Then,
\begin{equation}
\label{spring}
\M\big(\Delta(\mathbb{B}^n(a),\mathbb{S}^{n-1}(b);\mathbb{B}^n(b))\big)
= \omega_{n-1}\Big(\log \frac{b}{a}\Big)^{1-n}.
\end{equation}

\subsection*{Canonical ring domains}

The complementary components of the {\it Gr\"otzsch ring}
$R_{G,n}(s)$ in $\overline{\mathbb{R}}^n$ are
$\overline{\mathbb{B}}^n$ and $[se_1,\infty]$, $s>1$, and those of
the {\it Teichm\"uller ring} $R_{T,n}(s)$ are $[-e_1,0]$ and
$[se_1,\infty]$, $s>0$. We define two special functions
$\gamma_n(s)$, $s>1$, and $\tau_n(s)$, $s>0$, by
$$
\left\{ \begin{array} {llllll}
\gamma_n(s)=\M\big(\Delta(\overline{\mathbb{B}}^n,[se_1,\infty])\big)=\gamma(s),\\
\tau_n(s)=\M\big(\Delta([-e_1,0],[se_1,\infty])\big)=\tau(s),
\end{array} \right.
$$
respectively. The subscript $n$ is omitted if there is no danger of
confusion. We shall refer to these functions as the {\it Gr\"otzsch
capacity} and the {\it Teichm\"uller capacity}. It is well-known
that for all $s>1$
$$
\gamma_n(s)=2^{n-1}\tau_n(s^2-1),
$$
and that $\tau_n\colon(0,\infty)\to(0,\infty)$ is a decreasing
homeomorphism. For $s>1$ we have the following inequalities (see
e.g. \cite[7.24]{Vuorinen:1988}):
\begin{equation}
\label{qr7.24} \omega_{n-1}\big(\log\lambda_n s\big)^{1-n}\leq
\gamma(s) \leq \omega_{n-1}\big(\log s\big)^{1-n},
\end{equation}
where $\lambda_n$ is the Gr\"otzsch ring constant depending only on
$n$. The value of $\lambda_n$ is known only for $n=2$, namely
$\lambda_2=4$. For $n\geq 3$ it is known that $2e^{0.76(n-1)} <
\lambda_n \leq 2e^{n-1}$. For more information on the constant
$\lambda_n$, see \cite[Chapter 12]{Anderson:1997}.

We will use the following estimate from \cite{Gehring} (see also \cite[2.11]{MRV}).
Suppose that $G=A\setminus C$ is a ring domain such that $A\subset
\mathbb{B}^n$ and $C$
 is a
connected set with $0,x\in C$. Then
\begin{equation}
\label{mrv2.11} \M(\Delta(C, \partial A; G)) \geq \gamma(1/|x|).
\end{equation}

\subsection*{$K_I$- and $K_O$-inequalities}

Next we give two very useful inequalities, known as $K_I$- and
$K_O$-inequalities, respectively. The $K_I$-inequality is also known
as V\"ais\"al\"a's inequality.

\begin{thm}
$($\cite[Theorem II.9.1]{Rickman:1993}$)$ \label{vai_ie} Let
$f\colon G\to\mathbb{R}^n$ be a nonconstant quasiregular mapping,
$\Gamma$ be a path family in $G$, $\Gamma'$ be a path family in
$\mathbb{R}^n$, and $m$ be a positive integer such that the
following is true. For every path $\beta\colon I\to \mathbb{R}^n$ in
$\Gamma'$ there are paths $\alpha_1,\ldots,\alpha_m$ in $\Gamma$
such that $f\circ\alpha_j\subset\beta$ for all $j$ and such that for
every $x\in G$ and $t\in I$ the equality $\alpha_j(t)=x$ holds for
at most $i(x,f)$ indices $j$. Then
$$
\M(\Gamma')\leq \frac{K_I(f)}{m}\M(\Gamma).
$$
\end{thm}

In particular, we have the  Polecki{\u\i} inequality:

\begin{thm}
\label{poletskii} $($\cite[Theorem 8.1]{Rickman:1993}$)$ Let
$f\colon G\to\mathbb{R}^n$ be a nonconstant
 quasiregular mapping and let $\Gamma$ be a path family in $G$. Then
\[
\M(f\Gamma) \leq K_I(f)\M(\Gamma).
\]
\end{thm}

\begin{thm}
$($\cite[Theorem II.2.4]{Rickman:1993}$)$ \label{qr_thm1013} Let
$f\colon G\to \mathbb{R}^n$ be a nonconstant $K$-quasiregular
mapping. Let $A\subset G$ be a Borel set with $N(f,A)<\infty$, and
let $\Gamma$ be a family of paths in $A$. Then
$$
\M(\Gamma)\leq K_O(f)N(f,A)\M(f\Gamma).
$$
\end{thm}

\section{Topological properties}

Next we recall some topological properties of quasiregular mappings.

\subsection*{Discrete and open mappings.}

It is well-known that a nonconstant quasiregular mapping is discrete
 and
open. We denote by $B_f$ the branch set of $f$, i.e. the set of
points where $f$ fails to be a local homeomorphism. A result by V.\
A.~Chernavskii states that $\dim B_f \leq n-2$ for a discrete and
open $f\colon G\to\mathbb{R}^n$. The properties of discrete and open
mappings were further studied in by J.~V\"ais\"al\"a in
\cite{Vaisala:1966}, where also the multiplicity of discrete, open
and closed mappings was
 studied.

\subsection*{Normal domains.}

Let $f\colon G\to \mathbb{R}^n$ be a discrete and open mapping. A
domain $D\subset\subset G$ is called a {\it normal domain} for $f$ if
$f\partial D=\partial f D$. A normal neighborhood of $x$ is a normal
domain $D$ such that $D\cap f^{-1}(f(x))=\{x\}$.

\subsection*{Multiplicity and normal domains}

Let $f\colon G\to\mathbb{R}^n$ be a discrete and open mapping. We
denote 
\[
i(x,f) = {\inf_{U}}\, {\sup_{y}}\, \mathrm{card}\,f^{-1}(y)\cap U,
\]
where $U$ runs through the neighborhoods of $x$. The number
$i(x,f)$ is called the local (topological) index of $f$ at $x$. Let
$C\subset G$.
 The minimal multiplicity $M(f,C)$ and the maximal multiplicity
 $N(f,C)$
 are defined by
\begin{eqnarray}
M(f,C) &=& \inf_{y\in fC} \sum_{x\in f^{-1}(y)\cap C} i(x,f),\\
N(f,C) &=& \sup_{y\in fC} \sum_{x\in f^{-1}(y)\cap C} i(x,f),
\end{eqnarray}
respectively.

The following result holds for discrete, open and sense-preserving
mappings:

\begin{lem}
$($\cite[Corollary II.3.4]{Rickman:1993}$)$ \label{ric_34} Let $f\colon G \to \mathbb{R}^n$
be discrete, open and sense-preserving, $D\subset G$ a normal domain for $f$,
$\beta\colon[a,b) \to fD$ a path and $m=N(f,D)$. Then there exist
paths $\alpha_j\colon[a,b)\to D$, $1\leq j\leq m$, such that
\begin{itemize}
\item[(1)]  $f\circ\alpha_j=\beta$,
\item[(2)]  $\mathrm{card}\{j:\alpha_j(t)=x\}=i(x,f)$ for $x\in D\cap
 f^{-1}\beta(t)$,
\item[(3)]  $|\alpha_1|\cup\ldots\cup|\alpha_m|
=D\cup f^{-1}|\beta|$,
\end{itemize}
where $|\alpha|$ stands for the locus of $\alpha$, i.e. the image
set $\alpha [a,b)$, and $a\leq t < b$.
\end{lem}

\subsection*{Cluster sets}

The cluster set of $f\colon G\to\mathbb{R}^n$ at a point
$b\in\partial G$ is the set $C(f,b)$ of all points $z\in
\overline{\mathbb{R}}^n$ for which there exists a sequence $(b_k)$
in $G$ such that $b_k\to b$, and $f(b_k)\to z$. Let
$$
C(f,E)=\bigcup_{b\in E}C(f,b)
$$
for a non-empty set $E\subset\partial G$, and $C(f)=C(f,\partial
G)$. A mapping $f$ is {\it closed} if $fA$ is closed in $fG$
whenever $A$ is closed in $G$ and {\it proper} if $f^{-1}Q$ is compact in
$G$, where $Q$ is compact in $fG$. If $C(f)\subset\partial fG$, $f$
is said to be {\it boundary-preserving}.

\subsection*{Discrete, open and closed mappings}

Next we recall some useful topological results for discrete, open, and
closed mappings.

\begin{thm}
\label{vu1_thm33} $($See \cite[5.5]{Vaisala:1966},
\cite[3.3]{Martio-Srebro:1975} and \cite[3.2--3.3]{Vuorinen:1976}$)$
Let $f\colon G\to \mathbb{R}^n$ be discrete and open. Then the
following conditions are equivalent:
\begin{itemize}
\item[(1)] $f$ is proper.
\item[(2)] $f$ is closed.
\item[(3)] $f$ is boundary-preserving.
\item[(4)] Each sequence of points of $G$ converging to a point of
 $\partial G$ is transformed by $f$ onto a sequence no subsequence of
 which
 converges to a point of $fG$.
\item[(5)] $N(f,G)=p<\infty$ and for all $y\in fG$, we have
$$
p=\sum_{j=1}^k i(x_j,f),\,\,\, \{x_1,\ldots,x_k\}=f^{-1}(y).
$$
\end{itemize}
\end{thm}

\begin{cor}
\label{boundpre} If $f\colon G\to\mathbb{R}^n$ is discrete, open,
and closed, then $C(f)=\partial fG$.
\end{cor}

\begin{lem}
\cite[Lemma 3.6]{Vuorinen:1976} \label{vu1_lemma36} Let
$f\colon G\to \mathbb{R}^n$ be discrete, open, and closed, let
$U\subset fG$ be a domain, and let $D$ be a component of $f^{-1}U$.
Then $fD=U$ and $f|D$ is closed. Moreover, $C(f|D)=\partial U$. If
$f$ has a continuous extension $\overline{f}$ to $\overline{D}$,
then $\overline{f}\partial D=\partial U$.
\end{lem}

\begin{rem}
\label{cl-analytic}
In the plane each closed quasiregular mapping $f\colon \mathbb{B}^2\to \mathbb{B}^2$ has a representation
\[
f = g\circ h,
\]
where $h\colon \mathbb{B}^2\to \mathbb{B}^2$ is a quasiconformal mapping and $g\colon \mathbb{B}^2\to \mathbb{B}^n$ is a finite Blaschke product or a constant (see \cite[Theorem 4.1]{Vuorinen:1976}). This result follows immediately from the Sto\"ilov decomposition and the fact that each closed analytic function is a finite Blaschke product.
\end{rem}

\section{Unions of balls}

In this section, we prove a result which shows that a domain which is a union of a finite number of balls is always a $K$-quasiconformal image of a ball. The proof of this result also gives an explicit upper bound for the dilatation $K$.

We say that a domain $G\subset \overline{\mathbb{R}}^{n}$ is a {\it $K$-quasiball}, or simply {\it quasiball}, if there exists a $K$-quasiconformal mapping $f$ of $\overline{\mathbb{R}}^{n}$ onto itself such that
$G=f(\mathbb{B}^{n})$, where $\overline{\mathbb{R}}^n=\mathbb{R}^n\cup \{\infty\}$.

\begin{thm}\label{thm1.12}
Let $B_1,B_2,\ldots,B_m$ be balls in $\mathbb{R}^n$ such that for
 $1\leq j <m$, $|r_{j+1}-r_j|<|x_{j+1}-x_j|<r_j+r_{j+1}$ and $\overline{B_j}\cap
 \overline{B_k}=\emptyset$ for $|j-k|>1$. Then $D=B_1\cup B_2\cup \ldots \cup
 B_m$ is a quasiball.
\end{thm}

\subsection*{Wedge-shaped domains}

Let $(r,\varphi,z)$ be the cylindrical coordinates of a point $x\in
 \mathbb{R}^n$, $n\geq 3$. For $r\geq 0$,
 $0\leq\varphi < 2\pi$ $($ or $ -\pi \leq\varphi < \pi)$ and
$z\in \mathbb{R}^{n-2}=\{(0,0,z_3,\cdots, z_n):\; z_i\in
\mathbb{R},\; i=3,\cdots, n\}$ we define
$$
\left\{\begin{array}{rcl}
x_1 &=& r\cos \varphi,\\
x_2 &=& r\sin \varphi,\\
x_i &=& z_{i}\;\text{ for }\; 3 \leq i \leq n.
\end{array}\right.
$$
The domain $W_{(\gamma, \,\gamma+\alpha)}$, defined by
$\gamma<\varphi < \gamma+\alpha$, is called a {\it
 wedge} of angle $\alpha$, where
 $0\leq\gamma < 2\pi$, $0<\alpha < 2\pi$ and $0<\gamma+\alpha \leq 2\pi$
 (or $-\pi\leq\gamma < \pi$, $0<\alpha < 2\pi$ and $-\pi<\gamma+\alpha \leq \pi$).
 We also say that the domain
 $W_{(\gamma, \,\gamma+\alpha)}$ is a {\it
 wedge} of angle $\alpha$ with the starting angle $\gamma$.
For any rotation $\sigma$ around the subspace
 $\mathbb{R}^{n-2}$, $\sigma (W_{\gamma,\, \gamma+\alpha})$ is still a wedge of
angle $\alpha$. In particular,
 $W_{(\gamma,\, \gamma+\pi)}$ is a half-space in $\mathbb{R}^n$ for any $\gamma$.

Given two wedges $W_{(\gamma_1,\,\gamma_1+\alpha)}$ and
$W_{(\gamma_2,\,\gamma_2+\beta)}$, the quasiconformal
 diffeomorphism $f\colon W_{(\gamma_1,\,\gamma_1+\alpha)}\to W_{(\gamma_2,\,\gamma_2+\beta)}$, defined by
\[
f(r,\varphi,z)=(r,\beta\varphi/\alpha+\gamma_2-\gamma_1,z),
\]
is called a {\it folding}. Assuming that $0<\alpha\leq \beta<2\pi$ we have
$$
K_I(f) = \beta/\alpha, \qquad K_O(f)=(\beta/\alpha)^{n-1}.
$$
Then $f$ is a $(\beta/\alpha)^{n-1}$-quasiconformal mapping. See
\cite[16.3]{Vaisala:1971} for more details.

In what follows, we always denote by $\mathbb{B}^n(x_i, r_i)$ the
ball in $\mathbb{R}^n$ with the center $x_i$ and the radius $r_i$.

\begin{lem}\label{tme5.1}
\label{wedge} Suppose that $B_1$ and
 $B_2$ are two balls which satisfy $|r_2-r_1|<|x_1-x_2|<r_1+r_2$
 in $\mathbb{R}^n$.
 Then there exists $\alpha\in(\pi,2\pi)$ such that the domain $D=B_1\cup B_2$
 can be mapped onto a wedge $W_{(\gamma,\,\gamma+\alpha)}$ by a M\"obius transformation.
\end{lem}

\proof  Choose three distinct points $y_1,y_2,y_3\in S=\partial B_1\cap \partial B_2$. Then
 there exists (see e.g. \cite[7.21]{Anderson:1997}) a M\"obius transformation $g$
 such that $g(y_1)=0$, $g(y_2)=e_n$ and $g(y_3)=\infty$. It follows
 that $H_1=g(B_1)$ and $H_2=g(B_2)$ are half spaces in $\mathbb{R}^n$ and
 $0\in S'= \partial H_1 \cap \partial H_2$. We may assume that $S'$ is orthogonal to
 the $x_1x_2$-plane. Clearly $g(D)$ is a
 wedge $W_{(\gamma,\,\gamma+\alpha)}$ for some $\alpha \in (\pi,2\pi)$.
\qed

\subsection*{Angle of intersection }

Suppose that $B_1,B_2$ are two balls in $\mathbb{R}^n$ with
$|r_2-r_1|<|x_1-x_2|<r_1+r_2$. Then the angle of intersection,
$\alpha(B_1,B_2)$, of
 $B_1$ and $B_2$ is the number $\alpha \in (\pi,2\pi)$ such that there
 exists a M\"obius transformation $g$ such that $g(B_1\cup B_2)$ is a
 wedge $W_{(\gamma,\,\gamma+\alpha)}$ of angle $\alpha$.

\begin{cor}
Suppose that $B_1,\; B_2$ are balls in $\mathbb{R}^n$ with
$|r_2-r_1|<|x_1-x_2|<r_1+r_2$. Then $D=B_1\cup B_2$ is a
$K$-quasiball, where $K<\infty$ is
 a constant depending only on $\alpha(B_1,B_2)$ and $n$.
\end{cor}

\proof By Lemma \ref{wedge}, it is sufficient to prove that for any
$\alpha\in (\pi, 2\pi)$, the wedge $W_{(\gamma,\,\gamma+\alpha)}$ of
angle $\alpha$ is a
 quasiball.
Without loss of generality, we may assume that
$W_{(\gamma,\,\gamma+\alpha)}=W_{(0,\,\alpha)}$. Then the interior
of $\mathbb{R}^n\backslash W_{(0,\,\alpha)}$ is the wedge
$W_{(\alpha,\,2\pi)}$.
Let
$$
f(r,\varphi,z)= \left\{\begin{array}{rl}
(r,\pi\varphi/\alpha,z) & \text{ for } 0\leq \varphi \leq \alpha,\\
(r,\pi(1+(\varphi-\alpha)/(2\pi-\alpha)),z) & \text{ for } \alpha <
 \varphi < 2\pi,
\end{array}\right.
$$
and let $f(\infty)=\infty$. Then $f$ is clearly a homeomorphism of $\overline{\mathbb{R}}^n$
onto itself. It follows that $f\colon
 \overline{\mathbb{R}}^n\to\overline{\mathbb{R}}^n$ is quasiconformal with
 $K(f) = \max \{K(f_1),K(f_2)\}$, where the mappings $f_1\colon W_{(0,\, \alpha)}\to W_\pi=f_1(W_{(0,\, \alpha)})$ and $f_2\colon W_{(\alpha,\, 2\pi)} \to W_{(\pi,\, 2\pi)}=f_2(W_{(\alpha,\, 2\pi)})$, are foldings.
\qed

\proof[Proof of Theorem \ref{thm1.12}]
By Lemma \ref{wedge} we may find a M\"obius transformation $g$
taking $B_1\cup B_2$ onto a wedge $W_{(\gamma,\, \gamma+\alpha_1)}$
of angle $\alpha_1$ for some $\alpha_1\in (\pi,2\pi)$. Further, we
assume that $g(B_2)=W_{(0,\, \pi)}$ and $g(B_1)=W_{(\pi-\alpha_1,\,
2\pi-\alpha_1)}$, i.e., $W_{(\gamma,\,
\gamma+\alpha_1)}=W_{(\pi-\alpha_1,\, \pi)}$.
 Let $D_{m-i+1}=B_i\cup B_{i+1}\cup \cdots \cup B_m$ and
 $D'_{m-i+1}=g(D_{m-i+1})$ $(i=1, 2, \cdots, m)$. Define
$$
\varphi_0=\min\{ \varphi:\;  D'_{m-2}\backslash
\overline{W}_{(\pi-\alpha_1,\, \pi)}\subset W_{(\pi,\,
\pi+\varphi)}\}.
$$
Obviously, $0< \varphi_0<2\pi-\alpha_1$. We define a function
 $f_0\colon \overline{\mathbb{R}}^n \to \overline{\mathbb{R}}^n$ by
$$
f_0(r,\varphi,z)= \left\{\begin{array}{rl}
(r,\frac{\pi}{\alpha_1}(\varphi+(\alpha_1-\pi)),z) & \text{ for } \pi-\alpha_1<\varphi \leq \pi,\\
(r,\varphi,z) & \text{ for
 } \pi < \varphi \leq \pi+\varphi_0,\\
(r,\frac{(\pi-\varphi_0)\varphi+(\pi+\varphi_0)(\pi-\alpha_1)}
{(\pi-\alpha_1)+(\pi-\varphi_0)},z) & \text{ for } \pi+\varphi_0 <
\varphi \leq 3\pi-\alpha_1.
\end{array}\right.
$$
Then $f_0$ is a $K_1$-quasiconformal mapping, and $K_1<\infty$
depends
 only on $\alpha_1,\varphi_0$ and $n$. Let $f_1=g^{-1}\circ f_0 \circ
 g$. Then $f_1\colon \overline{\mathbb{R}}^n\to \overline{\mathbb{R}}^n$
 is $K_1$-quasiconformal and $f_1(D_m)=D_{m-1}$.

Similarly, for $j=2,\dots, m-1$ we may define a
$K_{j}$-quasiconformal
 mapping $f_j\colon \overline{\mathbb{R}}^n\to \overline{\mathbb{R}}^n$
 with $f_j(D_{m-j+1})=D_{m-j}$. Then
 $f=f_{m-1}\circ f_{m-2}\circ \cdots \circ f_1\circ h$ is a $K$-quasiconformal mapping
 of the whole space onto itself and $f(D_m)=\mathbb{B}^n$, where $h$
 is a suitable M\"obius transformation and
 $K=\prod_{j=1}^{m-1}K_j$. The claim follows.\qed

\section{Closed quasiregular mappings}

In this section, we study some of boundary regularity conditions, introduced by J.~V\"ais\"al\"a, under closed quasiregular mappings. These conditions are closely related to the boundary the mapping problems. We show that under certain assumptions boundary regularity conditions are preserved under closed quasiregular mappings. Indeed, without additional assumptions, the mapping properties of quasiregular mappings can be very different from quasiconformal ones, as illustrated by the following simple example.

\beg
It is well-known that one may map the unit ball $\mathbb{B}^n$ quasiconformally onto the half-ball $\mathbb{B}^n_+ = \{ x : |x|<1\textrm{ and } x_1>0\}$. Denote by $f_1$ a quasiconformal mapping such that $f_1\colon \mathbb{B}^n\to \mathbb{B}^n_+=f_1(\mathbb{B}^n)$. Let $f_2$ be the winding mapping defined by
\[
f_2(x_1,\ldots,x_n)= \big(r_2\cos(2\varphi_2),r_2\sin(2\varphi_2),x_3,\ldots,x_n\big),
\]
defined in the cylindrical coordinates such that $r_2=\sqrt{x_1^2+x_2^2}$ and $\varphi_2 = \arctan(x_2/x_1)$ for $x\in\mathbb{R}^{n}$. The mapping $f_2\colon \mathbb{R}^n\to \mathbb{R}^n$ is a well-known example of a quasiregular mapping (see e.g. \cite[3.1]{Rickman:1993}). In particular, we have $f_2(\mathbb{B}^n_+)= \mathbb{B}^n\setminus \{x: x_1\le 0,\,x_2=0\}$. Similarly, let $f_3\colon \mathbb{R}^n\to \mathbb{R}^n$ be the winding mapping defined by
\[
f_3(x_1,\ldots,x_{n}) = \big(r_3\cos(2\varphi_3),x_2,r_3\sin(2\varphi_3),x_4,\ldots,x_{n}\big),
\]
where $r_3=\sqrt{x_1^2+x_3^2}$ and $\varphi_3 = \arctan(x_3/x_1)$ for $x\in\mathbb{R}^n$. Then the quasiregular mapping $f=f_3\circ f_2 \circ f_1$ maps the unit ball onto the domain $\mathbb{B}^n\setminus \{x\in \mathbb{R}^n : x_1\le 0,\,x_2=x_3=0\}$.

\begin{center}
\includegraphics[width=11cm]{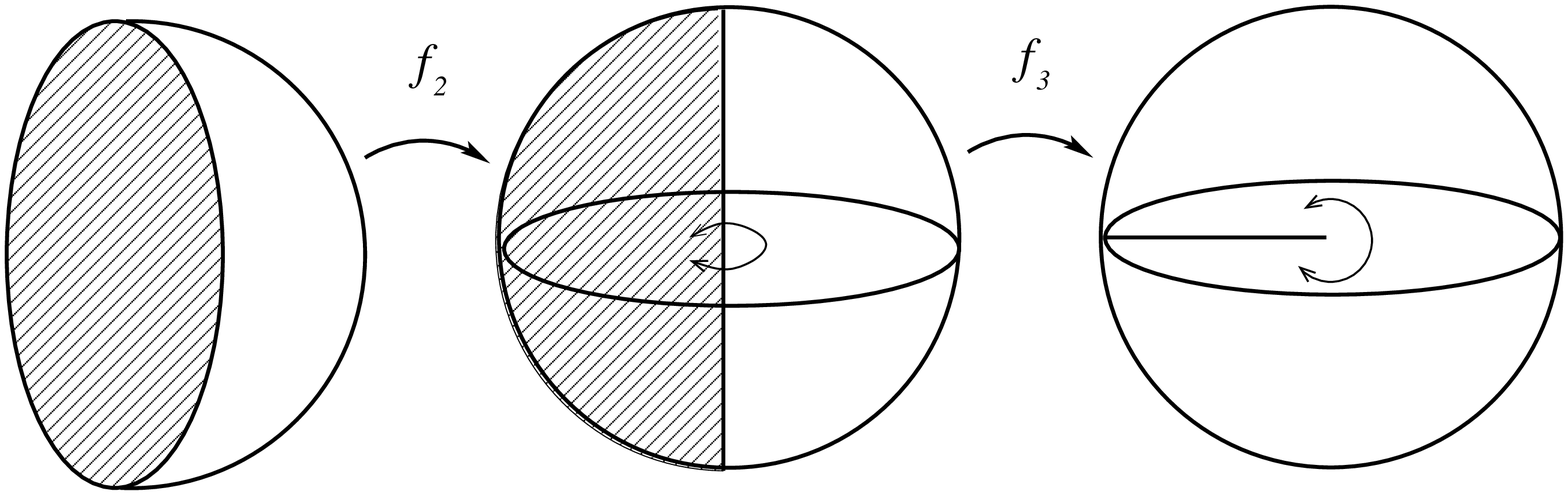}
\end{center}

In particular, for $n=3$ the image set is the unit ball with the negative $x_1$-axis removed. However, the cluster set of this mapping clearly consists of the unit sphere ${\mathbb S}^{2}$ and the two-dimensional disk
\[
\mathbb{D}=\{x\in \mathbb{R}^3: \sqrt{x_1^2+x_3^2}\le 1 \textrm{ and }x_2=0\},
\]
and thus the mapping $f$ is not closed.
\eeg

Our results in this section, Theorems \ref{thm3.1} and \ref{thm3.2}, are generalizations of similar results for
quasiconformal mappings (see \cite{Vaisala:1971}).

\subsection*{Boundary regularity conditions}

Recall that a quasiconformal map of $\mathbb{B}^n$ onto $\mathbb{B}^n$ has a
homeomorphic extension to $\overline{\mathbb{B}}^n$, see
\cite[Theorem 2]{Vaisala:1961}. The following definition is from \cite[17.5]{Vaisala:1971}.

\begin{defin}
Let $G$ be a domain in $\overline{\mathbb{R}}^n$ and let $b\in\partial G$.
\begin{itemize}
\item[(1)]
The domain $G$ is {\it locally connected} at $b$ if $b$ has
arbitrarily small neighborhoods $U$ such that $U\cap G$ is connected.
\item[(2)]
The domain $G$ is {\it finitely connected} at $b$ if $b$ has
 arbitrarily small neighborhoods $U$ such that $U\cap G$ has a finite
 number of
 components.
\item[(3)]
The domain $G$ has {\it property $P_1$} at $b$ if the following
 condition is satisfied: Whenever $E$ and $F$ are connected subsets of
 $G$ such
 that $b\in\overline{E}\cap\overline{F}$ we have
 $\M(\Delta(E,F;G))=\infty$.
\item[(4)]
The domain $G$ has {\it property $P_2$} at $b$ if: For each point
 $b_1\in\partial G$, $b_1\neq b$, there is a compact set
$F\subset G$, and a constant $\delta>0$, such that $\M(\Delta(E,F;
G))\geq\delta$ whenever $E$ is a connected set in $G$
 such that $\overline{E}$ contains $b$ and $b_1$.
\item[(5)] The domain $G$ is {\it locally quasiconformally collared} at
 $b$ if there is a neighborhood $U$ of $b$ and a homeomorphism $g$ of
 $U\cap \overline{G}$ onto the set $\{ x\in \overline{\mathbb{R}}^n :
 |x|<1\text{ and
 }x_n\geq 0\}$ such that $g|U \cap G$ is quasiconformal.
\item[(6)]
The domain $G$ is said to have one of the above properties at the
 boundary
if it has it at every boundary point.
\end{itemize}
\end{defin}

\begin{thm}\label{thm3.1}
 Suppose that $G$ and $G'$ are domains in
$\overline{\mathbb{R}}^n$, and let $f\colon \overline{G}\to
\overline{G'}$ be a continuous function such that $fG=G'$, and the
mapping $f_1=f|G$ is quasiregular and closed. If $G$ is a $P_1$
domain, and $G'$ is locally connected on the boundary, then $G'$ is
$P_1$.
\end{thm}

\begin{thm}\label{thm3.2}
 Suppose that $G$ and $G'$ are domains in
$\overline{\mathbb{R}}^n$, and let $f\colon \overline{G}\to
\overline{G'}$ be a continuous function such that $fG=G'$ and the
mapping $f_1=f|G$ is quasiregular and closed. If $G$ is a $P_2$
domain, then $G'$ also is $P_2$.
\end{thm}

Recall the next result from \cite[17.7]{Vaisala:1971}:

\begin{thm}\label{thm3.4}
 \label{vai71_177} The following
conditions are equivalent:
\begin{itemize}
\item[(1)] $G$ is finitely connected at $b$.
\item[(2)] Every neighborhood $U$ of $b$ contains a neighborhood $V$ of
 $b$ such that $V\cap G$ is contained in the union of a finite number
 of components of $U\cap G$.
\item[(3)] If $U$ is a neighborhood of $b$ and if $(x_j)$ is a sequence
 of points such that $x_j\to b$ and $x_j\in G$, then there is a
 subsequence which is contained in a single component of $U\cap G$.
\end{itemize}
\end{thm}

The next theorem, due to M.~Vuorinen, is a generalization of
\cite[17.13]{Vaisala:1971}.

\begin{thm}\label{thm3.5}
\cite[Theorem 4.2]{Vuorinen:1976} Suppose that $f\colon G\to
G'$ is a closed quasiregular mapping and that $G$ has the property
$P_1$ at the point $b\in\partial G$. Then the set $C(f,b)$ contains
at most one point at which $G'$ is finitely connected.
\end{thm}

The combination of Theorems \ref{thm3.4} and \ref{thm3.5} easily
implies the following result about the extension of closed
quasiregular mappings.

\begin{cor}
Let $f\colon G\to G'$ be a closed quasiregular mapping, let $G$ be a
$P_1$ domain, and let $G'$ be finitely connected on the boundary.
Then $f$ can be extended to a continuous mapping
$\overline{f}\colon\overline{G}\to\overline{G'}$.
\end{cor}

Now we are ready to prove Theorems \ref{thm3.1} and \ref{thm3.2}.

\proof[Proof of Theorem \ref{thm3.1}]
 Let $b'\in\partial G'$.
For a point $b$ in $\partial G$ we define a set $V(b,r)$ to be the
$b$-component of the set
$$
f^{-1}\big(\overline{G'}\cap\mathbb{B}^n(f(b),r)\big).
$$
Because $C(f_1)=\partial G'$ by Corollary \ref{boundpre}, we may
find a point $b\in\partial G$ such that $f(b)=b'$.

Now, let $E',F'$ be any continua in $G'$ such that
$b'\in\overline{E'}\cap\overline{F'}$. As $G'$ locally connected on
the boundary, each neighborhood $U$ of $b'$ is connected and
intersects with $E'$ and $F'$. Let $V$ be the $b$-component of
$f^{-1}U$. We choose $E,F$ to be the $b$-components of
$(f^{-1}E')\cap V\cap G$ and $(f^{-1}F')\cap V\cap G$, respectively.
It follows that $E,F$ are continua in $G$ and
$b\in\overline{E}\cap\overline{F}$. As $G$ is a $P_1$ domain,
$\M(\Delta(E,F;G))=\infty$. By Theorem \ref{vu1_thm33}(5),
$$
N(f,G)=p<\infty.
$$
Let $\Gamma=\Delta(E,F;G)$. By Theorem \ref{qr_thm1013}
$$
\M(\Gamma)\leq N(f_1,G)K_O(f_1)\M(f_1\Gamma),
$$
and thus
\begin{eqnarray*}
\infty & = &\M(\Gamma)              \\
&\leq & N(f,G)K_O(f_1)\M(f_1\Gamma) \\
&\leq & pK_O(f_1)\M(\Delta(E',F';G')).
\end{eqnarray*}
So, we have concluded that $\M(\Delta(E',F';G'))=\infty$, and the
claim is proved. \qed

\proof[Proof of Theorem \ref{thm3.2}] Let $b',
b_1'\in\partial G'$ such that $b_1'\neq b'$ and $E'\subset G'$ be a
connected set such that $b',b_1'\in\overline{E'}$.
By Lemma \ref{vu1_lemma36} we may choose
$$
E\subset(f^{-1} E')\cap G
$$
such that $fE=E'$ and $E$ is connected. As by Lemma
\ref{vu1_lemma36}
$$
f\partial E=\partial fE,\,\,\, f\partial G=\partial fG,
$$
and
$$
b',b_1'\in\partial E'\cap\partial G' =\partial fE\cap\partial fG.
$$
Hence, we may conclude that $\partial G\cap\partial E$ contains at least
two separate points, $b\in f^{-1}(b')$ and $b_1\in f^{-1}(b_1')$.

Now $b,b_1\in \partial G$ are separate points and $E$ is a continuum
such that $b,b_1\in \overline{E}$. It was assumed, that $G$ is a
$P_2$ domain, and so there exists a compact set $F$ and a constant
$\delta>0$ such that $\M(\Delta(E,F;G))\geq\delta$. As $f_1$ is a
closed quasiregular mapping, by Theorem \ref{vu1_thm33},
$N(f_1,G)=p<\infty$. By Theorem $\ref{qr_thm1013}$
\begin{eqnarray*}
\delta&\leq&\M\big(\Delta(E,F;G)\big)\\
&\leq&N(f_1,G)K_O(f_1)\M\big(\Delta(E',f_1F;G')\big).
\end{eqnarray*}
We may choose $F'=f_1F$ and
$$
\delta'=\frac{\delta}{pK_O(f_1)}>0.
$$
It follows that
$$
\M(\Delta(E',F';G'))\geq\delta'>0.
$$
As $f_1F$ is a compact set, the set $G'$ is a $P_2$ domain with the
corresponding compact set $F'$ and the constant $\delta'$, proving the
claim. \qed

The following problem related to the branch set of a closed quasiregular mappings was given by M. Vuorinen in 1980's \cite[p. 193]{Vuorinen:1988}, and it is still open.

\begin{prob}
Let $f\colon \mathbb{B}^n\to f\mathbb{B}^n\subset \mathbb{B}^n$ be discrete, open and proper. Assume that $n\ge 3$ and $B_f$ is compact. Is $f$ one-to-one? The answer is yes, if $f\mathbb{B}^n=\mathbb{B}^n$.
\end{prob}

\begin{rem}
A mapping $f\colon G \to \mathbb{R}^n$ is called {\it harmonic} if all its coordinate functions $u_j\colon G \to\mathbb{R}$ satisfy the Laplace equation $\Delta u_j = 0$. In particular, analytic functions are harmonic.   Recall that each closed analytic function is a finite Blaschke product or a constant (see Remark \ref{cl-analytic}). The class of harmonic mappings has been extensively studied \cite{Duren:2004}, and certain topological properties of harmonic mappings have been considered in \cite{Lyzzaik:1992}. However, to our knowledge, the class of closed harmonic mappings has not been studied.
\end{rem}

\section{Boundary behavior}

In this section, we prove some boundary behavior results for closed quasiregular mappings.

\subsection*{Existence of arcwise limits.}

A classical theorem by P.~Koebe states that a conformal mapping of a
simply connected domain $G$ in the complex plane $\IC$ has arcwise
limits along all end-cuts of $G$. R.~N\"akki \cite{Nakki:1979} proved a
similar result for quasiconformal mappings in $\mathbb{R}^n$. We show
that this result holds for closed quasiregular mappings as well.

Let $G$ a domain $\mathbb{R}^n$. A point $b\in \partial G$ is called
{\it accessible} from $G$ if there is a closed Jordan arc $\gamma$
contained in $G$ except for one endpoint, $b$. Then $\gamma$ is called an
{\it end-cut} of $G$ from $b$. Suppose that $f$ is a mapping of $G$ into
$\overline{\mathbb{R}}^n$. The cluster set of $f$ at $b$ along an end-cut
$\gamma$ from $b$ is denoted by $C_\gamma(f,b)$. If $C_\gamma(f,b)=\{
b'\}$, then $b'$ is called an {\it arcwise limit} of $f$ at $b$.

\begin{defin}
The spherical (chordal) metric $q$ in
 $\overline{\mathbb{R}}^n$ is
defined by
$$
\left\{\begin{array}{rcl} q(x,y)&=&
\frac{|x-y|}{\sqrt{1+|x|^2}\;\sqrt{1+|y|^2}}, \text{ for }
 x\neq\infty\neq y\;,\\
q(x,\infty)&=& \frac{1}{\sqrt{1+|x|^2}}.
\end{array}\right.
$$
For a set $E$ in $\overline{\mathbb{R}}^n$, we denote by $q(E)$ the
diameter of $E$ with respect to the metric $q(x,y)$.
\end{defin}

\begin{lem}
\label{q-est} $($\cite{Nakki:1973}$)$ Let $G$ be a locally
quasiconformally collared domain and let $E,F$ be
nondegenerate continua in $G$. Then for each $r>0$ there exists
$\delta>0$ such that $\M(\Delta(E,F;D))\geq \delta$ whenever $q(E)\geq r$ and $q(F)\geq r$.
\end{lem}

\begin{thm}
Suppose that $G$ is domain in $\mathbb{R}^n$, $f\colon G\to G'=fG$
is a closed quasiregular mapping, and $G'$ is a locally quasiconformally collared domain.
Then $f$ has arcwise limits along all end-cuts of $G$.
\end{thm}

\proof Let $b\in \partial G$, and suppose that $\gamma$ is an
end-cut from the
 point $b$. Fix a continuum $C \subset G$. We choose a sequence of
 neighborhoods $U_k$ of $b$ such that $\bigcap_{k=1}^\infty U_k= \{b\}$
 and
 $\gamma_k= U_k\cap G \cap \gamma$ is connected for $k=1,2,\ldots$. Write $C'=fC$. By Theorem \ref{vu1_thm33}, $f^{-1}C'$ is compact, and by Lemma \ref{ric_34} every path in $\Delta (C',|f(\gamma_k)|; G')$ has a lifting in $G$ beginning at $|\gamma_k|$ and leading to $f^{-1} C'$. Denote by $\Gamma_k$ the family of these liftings.
Then
\[
\lim_{k\to\infty} \M(\Gamma_k) = 0,
\]
and $f(\Gamma_k) < \Delta(C',|f(\gamma_k)|; G')$. Hence, by Theorem \ref{poletskii}, we have
\[
\M(\Delta(C',|f(\gamma_k)|; G')) \leq \M(f(\Gamma_k)) \leq K_I(f)\M(\Gamma_k) \to 0
\]
as $k\to\infty$. Then it follows by Lemma \ref{q-est} that $\lim_{k\to\infty}
 q(|f(\gamma_k)|)=0$ and hence $f$ has a limit at $b$ along $\gamma$.
\qed

\subsection*{Relative size of preimages}\label{sec4}

By Theorem \ref{vu1_thm33}, a set $D$ has at most $p<\infty$
preimages under a closed quasiregular mapping. Next we give an upper bound for the diameter of a preimage in
terms of the diameter of another preimage, i.e., we will prove that
only the images of the sets of roughly similar size can coincide in
a closed quasiregular mapping. Our result reads as follows.

\begin{thm}\label{thm4.1}

Let $f\colon G\to \mathbb{R}^n$ be a closed $K$-quasiregular
mapping. Suppose that $0<t<1$, and $A_1,A_2\subset
\mathbb{B}^n(x,tr)$ are nondegenerate continua with $A_1\cap
A_2=\emptyset$ such that $fA_1=fA_2$ and
$\overline{\mathbb{B}}^n(x,r)\subset G$. Then there is a
homeomorphism $h\colon[0,\infty)\to [0,\infty)$ depending only on
$n,K,t$ and $N(f,\mathbb{B}^n(x,r))$ such that $d(A_1)\geq
h(d(A_2))$.
\end{thm}

Before the proof of Theorem \ref{thm4.1}, we introduce two lemmas.

\begin{lem}
\label{grsphe} \cite[Lemma 2.31.]{Heikkala:2002} Let
$0<r_0<1$. Then
$$
C(n,r_0)\M\big(\Delta(\mathbb{B}(r),\mathbb{S}^{n-1})\big)\leq
\gamma_n(1/r) \leq
\M\big(\Delta(\mathbb{B}(r),\mathbb{S}^{n-1})\big)$$ for $r_0>r>0$,
where
$$
C(n,r_0)=\bigg(1-\frac{\log\lambda_n}{\log{r_0}}\bigg)^{1-n}.
$$
\end{lem}

\begin{lem}
\label{qr_143}
\cite[1.43]{Vuorinen:1988} Let $0<s<1$. Then
for all $a,x,y\in\overline{\mathbb{B}}^n(s)$
$$
\frac{1-s^2}{(1+s^2)^2}|x-y| \leq |T_ax-T_ay|\leq
\frac{1}{1-s^2}|x-y|.
$$
\end{lem}

\proof[Proof of Theorem \ref{thm4.1}] Let
$p=N(f,\mathbb{B}^n(x,r))<\infty$. By replacing $f$ with the mapping
$f \circ g$, where $g\colon z\mapsto (z-x)/r$, if necessary, we may
assume that $\mathbb{B}^n(x,r)=\mathbb{B}^n$. We choose the points
$z_1,z_2\in\overline{A_1}$ and $y_1,y_2\in \overline{A_2}$ such that
$d(A_1)\leq 2|z_1-z_2|$ and $d(A_2)\leq 2|y_1-y_2|$,
respectively. Next we estimate the modulus of curve family
$\Delta(A_1,\mathbb{S}^{n-1})$ with the capacity of spherical
annulus (\ref{spring}), and then apply 
Theorem \ref{poletskii} to obtain the estimate:
\begin{eqnarray*}
\omega_{n-1}\bigg[\log\bigg(\frac{1}{2|T_{z_1}(z_2)|}\bigg)\bigg]^{1-n}
& \geq &
\M\big(\Delta(A_1,\mathbb{S}^{n-1})\big)  \\
& \geq &
\frac{\M\big(f(\Delta(A_1,\mathbb{S}^{n-1}))\big)}{K_I(f)} \\
& = &
\frac{\M\big(\Delta(fA_1,f\mathbb{S}^{n-1})\big)}{K_I(f)} \\
& = & \frac{\M\big(f(\Delta(A_2,\mathbb{S}^{n-1}))\big)}{K_I(f)}.
\end{eqnarray*}
Now we apply the $K_O$-inequality, and then estimate the modulus in
terms of the capacity of the Gr\"otzsch ring domain
\begin{eqnarray*}
\frac{\M\big(f(\Delta(A_2,\mathbb{S}^{n-1}))\big)}{K_I(f)} & \geq &
\frac{\M\big(\Delta(A_2,\mathbb{S}^{n-1})\big)}{pK_I(f)K_O(f)} \\
& \geq & \frac{\gamma\big(|T_{y_1}(y_2)|^{-1}\big)}{pK_I(f)K_O(f)}.
\end{eqnarray*}
By combining these estimates with Lemma \ref{grsphe} and
(\ref{qr7.24}) we obtain
\begin{eqnarray*}
\omega_{n-1}\bigg[\log\bigg(\frac{1}{2|T_{z_1}(z_2)|}\bigg)\bigg]^{1-n}
&\geq &
\frac{\gamma\big(|T_{y_1}(y_2)|^{-1}\big)}{pK_I(f)K_O(f)} \\
&\geq & \frac{C(n,t)\omega_{n-1}}{pK_I(f)K_O(f)} \bigg[\log
\Big(\frac{\lambda_n}{|T_{y_1}(y_2)|}\Big)\bigg]^{1-n}.
\end{eqnarray*}
We have $(2\lambda_n|T_{z_1}(z_2)|\big)^{C(K,n,p,t)}\geq
|T_{y_1}(y_2)|$, and by applying Lemma \ref{qr_143} we obtain
$$
\bigg[2\lambda_n\frac{1-t^2}{(1+t^2)^2}|z_1-z_2|\bigg]^{C(K,n,p,t)}
\geq \frac{1}{1-t^2}|y_1-y_2|,
$$
proving the claim. \qed

{\bf Acknowledgments.}
We are indebted to the anonymous referee for very valuable suggestions concerning the presentation of this paper.


\begin{thebibliography}{99}

\bibitem{Ahlfors:1981}
{\sc L.\ V.~Ahlfors:} {\it M\"obius transformations in several
dimensions}, {Lecture notes, University of Minnesota, Minneapolis,
MN}, 1981.

\bibitem{Beurling-Ahlfors}
{\sc A.~Beurling} and {\sc L.~Ahlfors}:
The boundary correspondence under quasiconformal mappings. 
{\it Acta Math.} {\bf 96} (1956), 125--142.

\bibitem{Anderson:1997}
{\sc G.\ D.~Anderson, M.\ K.~Vamanamurty} and {\sc M. Vuorinen:} {\it
Conformal invariants, inequalities and quasiconformal mappings},
{Wiley-Interscience}, 1997.

\bibitem{Duren:2004}
{\sc P.~Duren:}
{\it Harmonic Mappings in the Plane},
Cambridge University Press, 2004.

\bibitem{Gehring}
{\sc F.\ W.~Gehring:} {Symmetrization of Rings in Space},
{\it Trans. Amer. Math. Soc.} {\bf 101} (3) (1961), 499--519.

\bibitem{Gehring-Hag}
{\sc F.\ W.~Gehring} and {\sc K.~Hag}:
Reflections on reflections in quasidisks. {\it Papers on analysis}, 81--90, 
Rep. Univ. Jyv\"askyl\"a Dep. Math. Stat., 83, Univ. Jyv\"askyl\"a, Jyv\"askyl\"a, 2001. 

\bibitem{Gehring-Vaisala}
{\sc F.\ W.~Gehring} and {\sc J.~V\"ais\"al\"a}: The coefficients of quasiconformality of domains in space, {\it Acta Math.} {\bf 114} (1965) 1--70.

\bibitem{Hag-Vamanamurthy}
{\sc K.~Hag} and {\sc M.\ K.~Vamanamurthy}:
The coefficients of quasiconformality of cones in n-space. 
{\it Ann. Acad. Sci. Fenn. Ser. A I} {\bf 3} (2) (1977), 267--275. 

\bibitem{Heikkala:2002}
{\sc V. Heikkala:} {Inequalities for conformal capacity, modulus
and conformal invariants}, {\it Ann. Acad. Sci. Fenn. Math.
Dissertationes} {\bf 132} (2002), 1--62.

\bibitem{Lyzzaik:1992}
{\sc A.~Lyzzaik:}
Local properties of some light harmonic mappings, {\it Canad. J. Math.} {\bf 44} (1992), 135--153.

\bibitem{MRV}
{\sc O.~Martio, S.~Rickman} and {\sc J.~V\"ais\"al\"a:} {Distortion
and singularities of quasiregular mappings}, {\it Ann. Acad. Sci. Fenn. A I Math.} {\bf 465} (1970), 1--12.

\bibitem{Martio-Srebro:1975}
{\sc O.~Martio} and {\sc U.~Srebro:} { Periodic quasimeromorphic
mappings}, {\it J. Anal. Math.} {\bf 28} (1975), 20--40.

\bibitem{Nakki:1973}
{\sc R.~N\"akki:} {Extension of Loewner's capacity theorem},
{\it Trans. Amer. Math. Soc.}
 {\bf 180} (1973), 229--236.

\bibitem{Nakki:1979}
{\sc R.~N\"akki:} {Prime ends and quasiconformal mappings}, 
{\it J. Anal. Math.} {\bf 35} (1979), 13--40.

\bibitem{Rickman:1993}
{\sc S.~Rickman:} {\it Quasiregular Mappings}, {Ergeb. Math.
Grenzgeb. {\bf 3}, Vol. 26, Springer-Verlag}, {Berlin}, 1993.

\bibitem{Vaisala:1961}
{\sc J.~V\"ais\"al\"a:} {On quasiconformal mappings of a ball},
{\it Ann. Acad. Sci. Fenn. Ser. A I } {\bf 304} (1961), 1--7.

\bibitem{Vaisala:1966}
{\sc J.~V\"ais\"al\"a:} {Discrete and open mappings on
manifolds}, {\it Ann. Acad. Sci. Fenn. Ser. A I} {\bf 392} (1966),
1--10.

\bibitem{Vaisala:1971}
{\sc J.~V\"ais\"al\"a:} {\it Lectures on n-Dimensional
Quasiconformal Mappings}, {Lecture Notes in Math., Vol. 229,
Springer-Verlag}, {Berlin}, 1971.

\bibitem{Vuorinen:1976}
{\sc M.~Vuorinen:} {Exceptional sets and boundary behavior of
quasiregular mappings in $n$-space}, {\it Ann. Acad. Sci. Fenn. Ser. A I Math.
Dissertationes} {\bf 11} (1976),
 1--44.

\bibitem{Vuorinen:1978}
{\sc M.~Vuorinen:} {On angular limits of closed quasiregular
mappings}, {\it Proceedings of the {F}irst {F}innish-{P}olish {S}ummer
              {S}chool in {C}omplex {A}nalysis ({P}odlesice, 1977), {P}art
              {II}}, Univ. \L \'od\'z, \L \'od\'z,  1978, 69--74.

\bibitem{Vuorinen:1985}
{\sc M.~Vuorinen:} {Conformal invariants and quasiregular
mappings}, {\it J. Anal. Math} {\bf 45} (1985), 69--115.

\bibitem{Vuorinen:1988}
{\sc M.~Vuorinen:} {\it Conformal Geometry and Quasiregular
Mappings}, {Lecture Notes in Math., Vol. 1319, Springer-Verlag},
{Berlin}, 1988.


\end{thebibliography}
\end{document}